\documentclass[12pt,reqno]{amsart}

\usepackage{mathrsfs}
\usepackage{multicol}
\usepackage{fancyvrb}
\usepackage{multicol}
\usepackage{fancyvrb}
\usepackage{xypic}
\usepackage{amssymb,latexsym,amsbsy,amsthm,amsxtra,amscd,amsfonts,amsthm,amsmath,verbatim}
\usepackage{amsxtra,amsmath,epsfig}

\newcommand{\ncom}{\newcommand}
\ncom{\bq}{\begin{equation}}
\ncom{\eq}{\end{equation}}
\ncom{\beqn}{\begin{eqnarray*}}
\ncom{\eeqn}{\end{eqnarray*}}
\ncom{\beq}{\begin{eqnarray}}
\ncom{\eeq}{\end{eqnarray}}
\ncom{\been}{\begin{enumerate}}
\ncom{\eeen}{\end{enumerate}}
\ncom{\olin}{\overline}
\ncom{\f}{\frac}
\ncom{\rar}{\rightarrow}
\def\a{{\bf a}}

\usepackage{lmodern}
\usepackage[T1]{fontenc}

\usepackage{color}

\makeatletter

\@addtoreset{equation}{section}
\def\theequation{\thesection.\@arabic \c@equation}

\makeatother

\def\theequation{\arabic{equation}}
\def\theequation{\thesection.\arabic{equation}}
\numberwithin{equation}{section}

\def\charac{\operatorname{char}}

\def\dim{\operatorname{dim}}
\def\gcd{\operatorname{gcd}}

\def\height{\operatorname{ht}}

\def\ker{\operatorname{ker}}
\def\max{\operatorname{max}}

\def\reg{\operatorname{reg}}

\newcommand{\kk}{\Bbbk}

\def\A{\mathbb A}
\def\C{\mathcal C}

\def\Proj{\mathbb P}
\def\R{\mathcal R}

\def\a{\ulin{a}}

\newcommand{\p}{\mathfrak p}
\newcommand{\pp}{\mathtt p}

\newcommand{\m}{\mathfrak m}

\theoremstyle{plain}
\newtheorem{theorem}[equation]{Theorem}

\newtheorem{lemma}[equation]{Lemma}
\newtheorem{conjecture}[equation]{Conjecture}

\theoremstyle{definition}

\newtheorem{problem}[equation]{Problem}

\newtheorem{question}[equation]{Question}

\ncom{\bib}{\bibitem}

\ncom{\limns}{\underset{\underset{s}{\longrightarrow}}{\lim}}
\ncom{\limnr}{\underset{\underset{r}{\longrightarrow}}{\lim}}

\ncom{\maxi}{\underset{{i}}{\max}}

\ncom{\Tprime}{T^{\prime}}
\ncom{\mprime}{\m^{\prime}}

\def\ulin{\underline}

\begin{document}
 \title[ ] { Symbolic powers, set-theoretic complete intersection and certain invariants}
 \author{Clare D'Cruz}
 \address{Chennai Mathematical Institute, Plot H1 SIPCOT IT Park, Siruseri, 
Kelambakkam 603103, Tamil Nadu, 
India
} 
\thanks{$^*$ partially supported by a grant from Infosys Foundation}
\email{clare@cmi.ac.in}

\maketitle
\begin{abstract}
In this survey article we give a brief history of symbolic powers  and its  connection with the interesting problem of set-theoretic complete intersection. We also state a few  problems and conjectures. Recently, in connection to symbolic powers is the containment problem. We list a few interesting results and  related problems on the resurgence, Waldschmidt constant and Castelnuovo-Mumford regularity. 
\end{abstract}

\section{Introduction}
Let $A$ be a Noetherian ring and $I$ an ideal in $A$ with no embedded components. Then the ideal 
${\displaystyle I^{(n)} :=  A \cap (\bigcap_{\p \in Ass(A/I)}   I^n A_{\p} )}$ is the $n$-th symbolic power of $I$.  The study of $n$-th symbolic power has been important for the last few  decades mainly because  of its connection with algebraic geometry. In recent years, it has become even more active area of research mainly because several  interesting associated invariants. We list some of the interesting questions and open problems.

To understand the connection with algebraic geometry,   let $\kk$ be a field and let $\A^d$ (or  $\A_{\kk}^d$) denote the set or all $d$-tuples $\a= (a_1, \ldots, a_d)$ where $a_i \in \kk$ for all $i=1, \ldots, d$.  The set $\A^d$  is called the affine $d$-space of dimension $d$ over $\kk$. 
We say that  a subset   $Y$ in $\A^d$ is a zero  set if it is the set of common zeros of a collection of polynomials $f_1, \ldots, f_m \in R:= \kk[X_1, \ldots, X_d]$ and we denote it by $Y = Z(f_1, \ldots, f_m)$. We can define a topology on $\A^{n}$ by defining the closed sets to be the zero sets. If $I= (f_1, \ldots, f_m)$, then $Y = Z(I)$. To every subset of $Y \subset \A^n$ we can define the ideal of $I(Y):= \{  f \in \kk[X_1, \ldots, X_n] | f(P) = 0 \mbox{ for all } P \in Y\}$. An irreducible closed subset $Y$ of $\A^n$ called an affine algebraic set. 

Let $Y$ be an algebraic set.  We say that $Y$ is defined set-theoretically by $n$ elements if there exists  $n$ elements $f_1, \ldots, f_{n} \in R$ such that  $I(Y)= 
\sqrt{(f_1, \ldots, f_{n}})$. 
 In 1882, Kronecker showed that $I$ can be set theoretically defined by $d+1$ equations \cite{kronecker}. 
Later, this result was improved by  Storch \cite{storch} and by Eisenbud and Evans \cite{eisenbud-evans}. It follows from their work that if
$\kk$ is algebraically closed and $I$ is an homogenous ideal, then $I$ can be defined set-theoretically by $d$ 
elements.  Hence it was of interest to know which ideals $I$ could be defined set-theoretically by $d-1$ elements. 
 If $I$ is locally complete intersection  of pure dimension one,  then $I$ can be defined set-theoretically by  $d-1$ elements    (\cite{ferrand-1975},  \cite{bor-1978},  \cite{kumar-1978}, \cite{szpiro-1979}). 
In 1992,  Lyubeznick showed that    if $V$ is an algebraic set in $\A_{\kk}^d$ and $\charac(\kk)=\pp >0$, then   $V$ can be 
defined set-theoretically by $d-1$ elements \cite{lyubeznik-1992}.

In $1978$,  Cowsik and Nori  proved a remarkable result. They showed   that if $\charac(\kk) = \pp >0$, then  any affine curve is a set-theoretic complete intersection (\cite[Theorem~1]{cowsik-nori}).
If $\charac(\kk) =0$, then one of the best known results in $\charac(\kk)=0$ is the result of Herzog which was later also proved by Bresinsky \cite{bresinsky-1979}. He showed that all 
monomial curves in $\A^3$  are set-theoretic complete intersection.

In $1981$ Cowsik proved an interesting  result which connects commutative algebra and algebraic geometry: 
 \begin{theorem}
 \cite{cowsik-1981} Let 
$(R, \m)$ a Noetherian  local ring and $\p \not =\m $ a prime ideal. If  the symbolic Rees algebra 
$R_{s}(\p) := \oplus_{n \geq 0} \p^{(n)}$ is Noetherian, then $\p$
 can be defined set-theoretically by $d-1$ elements.
\end{theorem}
However, the converse need not be true. Cowsik's result motivated several researchers to investigate the Noetherian property of the symbolic Rees algebra. In 1987,  Huneke gave necessary and sufficient conditions for 
$\R_{s}(\p)$ to be Noetherian when $\dim~ R=3$ \cite{huneke-1987}. Huneke's result  was 
generalised in $1991$ for $\dim~ R \geq 3$ by Morales  \cite{morales-1991}. All these 
results paved a new  way to  study the famous problem on set theoretic complete 
intersection. In section~\ref{stci} we discuss some of  these problems.

We also have the famous result due to Zariski \cite{zariski} and Nagata  \cite{nagata} which states that  if $\kk$ is an algebraically closed  field, then the  $n$-th symbolic power of a given prime ideal consists of the elements that vanish up to order $n$ on the corresponding variety \cite{eisenbud}. This result has been generalised to perfect fields $\kk$ and radical ideals \cite[Proposition~2.14, Exercise~2.15]{dao}.

In general, symbolic powers of ideals are hard to compute. Hence recently, associated to symbolic powers of ideals  Bocci and Harbourne introduced a quantity called the resurgence \cite{BH0}. In section~\ref{resurgence}, we will state the recent developments on this quantity and related invariants like the Waldschmidt constant and Castelnuovo-Mumford regularity.

\section{Set-theoretic complete intersection and symbolic Rees algebra}
\label{stci}
\subsection{Set-theoretic complete intersection  in $\A^n$ and $\Proj^n$}\hfill\\
Throughout this section $R= \kk[X_1, \ldots,X_d]$ where $X_1, \ldots X_d$ are variables. 
We  say that a radical ideal $I \subset R$ is a set theoretic complete intersection of there exists  $h= \height(I)$ elements such that $I = \sqrt{(f_1, \ldots. f_h)}$.  
Let $\a:=(a_1, \ldots, a_d)$ be positive integers such that $\gcd(a_1, \ldots, a_d) =1$. Let  ${\C(\a) }:=\{ t^{\a} = 
(t^{a_1}, t^{a_2}, \ldots, t^{a_d} ) | t \in \kk\}$ be a curve in $\A^n$.  If $\phi: R \longrightarrow \kk[T]$ is 
the homomorphism given by $\phi(X_i) = T^{a_i}$ for all $i=1, \ldots, d$.   Then $ \p( \C(\a) ) := \ker(\phi)$ is the prime ideal defining 
the curve  $\C(\a)$. In other words, $I(\C(\a)) = \p( \C(\a) )$. We say that $\C(\a)$ can be defined set-theoretically by $d-1$ elements 
if there exists $d-1$ elements $f_1, \ldots, f_{d-1} \in \p$ such that $\p = \sqrt{(f_1, \ldots, f_{d-1})}$.

Let $d=3$, then  we have the interesting result: 
\begin{theorem}
Let $\gcd(a_1, a_2, a_3) =1$. 
\been
\item
 \cite{herzog-1970} Then  one of the following is true:
\been
\item

$\p( \C(\a) )$  is a complete intersection.
\item
There exists integers $\alpha_i, \beta_i, \gamma_i$; $i=1,2$ such that $\p(\C(\a))$ is  generated by $2 \times 2$ minors of the matrix  
${\displaystyle 
\left(\begin{array}{ccc} 
X_1^{\alpha_1} & X_2^{\beta_1} & X_3^{\gamma_1} \\
 X_2^{\beta_2} & X_3^{\gamma_2}  & X_1^{\alpha_2}  \\\end{array}\right) }$.
 \eeen
 \item
 {{\em [Herzog (Unpublised work)]} and \cite{bresinsky-1979}}
 $\p(\C(\a))$ is a set-theoretic complete intersection.
 \eeen
 \end{theorem}
 
Such a result is not known for $\A^d$, $d \geq 4$. The first result in higher dimension was given by Bresinsky where he showed that 
 certain Gorenstein curves in $\A^4$ are set-theoretic complete intersection \cite{bresinsky-a4-1979}. 
  In 1990,  Patil proved the following result \cite{patil-1990}:
  \begin{theorem}
  Let $a_1, a_2, \ldots, a_{d-2}$ be an arithmetic sequence. Then $\p(\C(\a))$ is a set-theoretic complete intersection.
  \end{theorem}
 In 1980, Valla showed that certain determinantal ideals were set-theoretic complete intersection \cite{valla-1980}. As a consequence they prove the following:
 
 \begin{theorem}
 \cite[Example~3.3]{valla-1980}
 Let $q,m$ be positive integers with $\gcd(q,m)=1$. Put $a_i = 2q+1 + (i-1)m$, for $i=1,2,3$. Then $\p(\C(\a))$ is a set-theoretic complete intersection. 
 \end{theorem}
 
 In the past forty years several researchers have given interesting examples of affine varieties which are set-theoretic complete intersection. However, the following question is still open. 

 \begin{question}
 Let $\kk$ be a field of characteristic zero and $d \geq 4$. Is every curve $\C(\a) \subseteq \A^d$ a set theoretic complete intersection?
 \end{question}
 
 We would like to bring to the attention a paper of Moh where he considered  the set-theoretic complete intersection problem of analytic space curves over an algebraically closed field. 
  Let $\kk[[X,Y,Z]]$ and $\kk[[T]]$ be power series rings and $\phi:\kk[[X,Y,Z]] \longrightarrow \kk[[T]]$ be given by 
 $\phi(X) = T^a + \cdots$,  $\phi(Y) = T^b + \cdots$ and  $\phi(Z) = T^c + \cdots$. Let $\p = \ker(\phi)$. Such curves are called Moh curves. Moh showed that if $(a-2)b <c $, then $\p$ is a set- theoretic complete intersection \cite{moh-1982}.   In \cite{huneke-1987} Huneke showed that the  symbolic Rees algebra of the Moh curve parameterized by $(t^6, t^7+t^{10}, t^8)$ is Noetherian. However, it is not easy to describe the defining ideal of a Moh curve. The following question is still open:
 
 \begin{question}
Let $\phi: \kk[[X,Y,Z]] \longrightarrow \kk[[T]]$ be given by 
 $\phi(X) = T^6 + T^{31} $,  $\phi(Y) = T^8 $ and  $\phi(Z) = T^{10}$. Let $\p = \ker(\phi)$. Is $\p$ a set-theoretic complete intersection?
 \end{question}

 We now focus our  attention on curves in the projective space $\Proj_k^n$. 
 It is a long standing question whether every connected subvariety in $\Proj_k^n$ is a set-theoretic complete intersection \cite{hartshorne-LNM-156}. The answer is not known even for curves in $\Proj^3$. We list a few results in this direction. 

Let $\a:=(a_1, \ldots, a_d)$ be integers such that $\gcd(a_1, \ldots, a_d) =1$ and 
$0= a_0<a_1 < a_2 \cdots < a_d$.  Put $S = \kk[X_0, X_1, \ldots, X_d] $ and let 
 $\psi: S \longrightarrow \kk[T,U]$ be
the homomorphism given by $\phi(X_i) = T^{a_d-a_i} U^{a_i}$ for all $i=0, \ldots, d$.  Then $ \p( {\olin{\C}}(\a) ) := \ker(\phi)$ is the prime ideal defining ideal of 
the curve  ${\olin{\C}}(\a)$. In other words, $I({\olin{\C}}(\a) = \p( {\olin{\C}}(\a) )$. 

One of the most simplest and interesting example was given by Hartshorne. 
\begin{theorem}
\cite{hartshorne-charp}
Let $\kk$ be a field of positive characteristic  $\pp$ and  $\a= (1, a-1,a)$, where $a \geq 4$. Then  $\p( {\olin{\C}}(\a) )$ is a set-theoretic complete intersection. 
\end{theorem}

Hartshorne's result was generalised by Ferrand \cite{ferrand-1975}. 
Robbiano and Valla studied the curve  ${\olin{\C}}(\a)$ for   $\a= (1, 3,4)$ \cite{robbiano-valla-1983}. 
In 1991, Moh generalised the work of Hartshorne and Ferrand \cite{moh-1985}. 

An interesting result in $\Proj^3$ is:

\begin{theorem}
\cite {robbiano-valla-1983}, \cite{stuckard-vogel-1982}
Let ${\mathcal C}$ be a curve in $\Proj^3$. Let $I({\mathcal C}) \subset S= \kk[X_0, X_1, X_2, X_3]$ be the ideal  of the ${\mathcal C}$. If $S/ I({\mathcal C})$ is Cohen-Macaulay, then $ I({\mathcal C})$ is a set-theoretic complete intersection.
\end{theorem}
 
A considerable amount of research has been done in this area. It is impossible to list all of them 
there.  The list of references is not exhaustive. For a good collection of articles on set-theoretic complete intersection 
one can read \cite{LNM-1092}. An interesting survey on this topic was also given by Lyubneznik 
\cite{lyubeznik-survey}. 

\subsection{Symbolic Rees Algebra} \hfill\\
The work of Cowsik, Huneke and Morales motivated several researchers to work on the 
 on the symbolic powers of a prime ideal and the symbolic Rees algebra $\R_{s}(\p) =\oplus_{n \geq 0 } \p^{(n)}$. 

If  $(R, \m)$ a Noetherian local ring of dimension $d$ and $\p$ a prime ideal of height $d-1$.  Some of the interesting questions on the symbolic Rees algebra are:
(1) Is it Noetherian? (2) Is it Cohen-Macaulay (3) Is it Gorenstein?
An answer to question (1) would imply that $\p$ is a set-theoretic complete 
intersection, by Cowsik's result.

 If $\phi: \kk[[X_1, \ldots, X_d]] \longrightarrow \kk[[T]]$ is 
the homomorphism given by $\phi(X_i) = T^{a_i}$ for all $i=1, \ldots, d$,  then $ \p(\a) := \ker(\phi)$.

In 1982, Huneke gave examples prime ideals $\p$ in $\kk[[X_1, X_2, X_3]]$ whose symbolic Rees algebra $R_s(\p)$ is Noetherian \cite{huneke-1982}. In fact he showed the following result.

\begin{theorem}
$R_s(\p(\a))$ is Noetherian in the following cases:
\been
\item
$\a = (2t+1, 2r+s, s+r+rs)$, where either  $s \leq r$ or $s >r$ and $t=1$.
\item
$\a = (s+2, 2r+1,  s+1+ rs)$, $2\leq r \leq s$.

\item
$\a= (t+s+1, tr+t+1,  rs+r+s)$, $r \leq t$ and $s\geq 1$.
\eeen
\end{theorem}

The first result which gave necessary and sufficient conditions for $R_s(\p)$ to be Noetherian  was by Huneke  \cite[Theorem~2.1]{huneke-1987}. This result was later generalised by Morales \cite[Theorem~2.1]{morales-1991}. A consequence of their result is:

\begin{theorem}
Let $(R, \m)$ be a regular local ring of dimension $d$ and $\p$ a prime ideal of height $d-1$. Then 
$R_s(\p)$ is Noetherian if and only if there exists $x \in \m \setminus \p$  and elements $f_i \in \p^{(k_i)}$, $i=1, \ldots, f_{d-1}$, 
such that 
\beqn
\ell \left( \f{R}{(\p, f_1, \ldots, f_{d-1},x)} \right)=\ell \left( \f{R}{ (p,x) }\right)   k_1 \cdots k_{d-1}.
\eeqn
\end{theorem}

Using this criteria several researchers have produced  examples of monomial curves $\p \in \kk[[X_1, X_2, X_3]]$ such that the  symbolic Rees algebra $R_s(\p)$ is Noetherian. We cite a few examples here \cite{schenzel-1988} , \cite{cutkosky-1991}, \cite{goto-nis-shim}
 \cite{hema-1991}, \cite{goto-nis-shim-2}, \cite{schenzel-1992}, \cite{goto},  \cite{sahin-2009-pn}.
 
One interesting question is: If $R$ is a Noetherian ring,  and $\p$ is a prime ideal, is the symbolic Rees algebra $R_s(\p)$ Noetherian? Rees provided an example which implies that this question does not have a positive answer \cite{rees-1958}. 
Later Cowsik conjectured that if $R$ is a regular local ring  and $\p$ is a prime ideal, then  $R_s(\p)$ is Noetherian. 
In 1990, Roberts gave a counter example  to Cowsik's conjecture \cite{roberts-1985}. In 1994, Goto, Nishida and Watanabe gave an infinite class of monomial curves whose symbolic Rees algebra  $R_s(\p)$ is Noetherian if the characteristic of $\kk$ is  $\pp$, but if the characteristic of $\kk$ is zero, then $R_s(\p)$ is not Noetherian \cite{goto-nis-wat}. 

We end this section stating a problem  which is still open:

\begin{problem}
 Let $\kk$ be an algebraically closed field.  Can one classify all monomial curves in  $\A^3_{\kk}$  for which $R_s(\p)$ is Noetherian?
\end{problem}

An interesting paper in this direction is \cite{cutkosky-1991}.

\section{Resurgence, Waldschmidt constant and regularity}
\label{resurgence}
One of the reasons the symbolic Rees algebra is hard to analyse is because it is not easy to describe the  symbolic powers even for curves in $\A^3_{\kk}$. Hence, one would like to compare the symbolic powers and ordinary powers of an ideal. If $I$ is an ideal in a Noetherian ring $R$, then   from the definition of symbolic powers it follows that  $I^n\subseteq I^{(n)}$ and in fact for any proper ideal nonzero ideal  $ I$,  $I^r\subseteq I^{(n)}$ holds if and only if $r \geq n$.  
A challenging problem to determine for which  positive integers $n$ and $r$ the containment  $I^{(n)}\subseteq I^r$ holds true. 
In \cite{swanson-2000}  Swanson compared the symbolic powers  and ordinary powers of several ideals. 

For the rest of this section we will assume that $S= \kk[X_0, \ldots, X_d]$ and $I$ is an homogenous ideal in $S$.
Hence in 2001, Ein, Lazarsfeld and Smith proved a very interesting result. It follows from their result  
\begin{theorem}
\cite{ein-laz-smith-2001}. Let $I \subset S$ be a proper ideal.  If $h$ is the largest height of an associated prime of $I$, then $I^{(hn)} \subseteq  I^n$ for all $n \geq  0$. 
\end{theorem}
In 2002, Hochster and Huneke proved a stronger result  \cite{ {hochster-huneke-2002}}. It follows from the above results that  if $d = \dim~R$, then  $I^{(n)}\subseteq I^r$ for $n\geq (d-1)r$. In this direction,  Harbourne raised the following: conjecture in:  

\begin{conjecture}
\cite[Conjecture~8.4.2]{BDHKKASS} For any homogeneous ideal $0 \not = I \subset S$, $I^{(n)} \subseteq I^r $ if $n \geq rd - (d-1)$. 
 \end{conjecture}
 In the same paper they remark that from the methods in \cite{hochster-huneke-2002} there is enough evidence for this conjecture to be true at least when the characteristic of $\kk$ is $\pp$ and $r = {\pp}^t$ for $t >0$ (see \cite[Example~8.4.4]{BDHKKASS}. 
 This has led to the study  of the least integer $n$ for which $I^{(n)}\subseteq I^r$ holds for a given ideal $I$ and  for an integer $r$. 
To answer this question C. Bocci and B. Harbourne  defined an asymptotic quantity \cite{BH-2010} called resurgence which is defined as
\beqn
\rho(I) := \sup\{ m/r \mid 
  I^{(m)}  \not \subseteq I^r\}.
  \eeqn
Hence if $m > \rho(I) r$, then $I^{(m)} \subseteq I^r$.  

In general resurgence is not easy to compute. Hence it is useful to  give bounds. 
From the results in \cite{ein-laz-smith-2001} it follows that $\rho(I) \leq d-1$.

Another interesting invariant is the Waldschmidt constant. 
This constant was introduced by Waldschmidt  in \cite{waldschmidt}. 
Let $I$ be an homogenous ideal and let   $ \alpha(I)$ denote  the least degree of a homogeneous generator of $I$. 
Then we have the famous conjecture due to Nagata:

\begin{conjecture}
\cite{nagata-1965} Let $V$ be a finite set of  $n$ points in $\Proj^{2}_{\C}$ and $I(V)$ be the corresponding homogenous ideal in $\C[X_0, X_1, X_2]$. Then 
$\alpha(I^{(m)}) \geq m \sqrt{n}$. 
\end{conjecture}
This conjecture is still open in general. It is know only in a few cases. 

Define
  $${\gamma(I) }:= \lim_{n\rightarrow \infty}\frac{\alpha(I^{(n)}) }{n}.$$
  $\gamma(I)$ is called the Waldschmidt constant (\cite{waldschmidt}, \cite{BH0}). Since $\alpha$ is subadditive, i.e., 
  $\alpha(I)^{(n + m)} \leq  \alpha(I^{(n)}) + \alpha(I^{(m)})$, it follows that $\gamma(I)$ exists (\cite{waldschmidt},  \cite{BH-2010}). Moreover,  there is a  lower bound for $\rho(I)$. It follows from 
    \cite[Lemma~2.3.2]{BH-2010}:
    
  \begin{lemma}
Let $0 \not = I \subset S$ be a homogenous ideal. Then $\gamma(I) \geq 1$ and 
\beqn
1 \leq \f{\alpha(I)}{\gamma(I)} \leq \rho(I).
\eeqn
  \end{lemma}
    Related to the Waldschmidt constant is the following conjecture:
 \begin{conjecture} 
 [Chudnovsky] 
 Let $V$ be a set of points in $\Proj^{d}$ and $I(V)$  be  the corresponding  homogenous ideal in $S$. Then
 \beqn
 \gamma(I) \geq \f{\alpha(I) + d-1}{d}. 
 \eeqn
  \end{conjecture}
  Recently,  Chudnovsky's conjecture has attracted the attention of researchers 
  (\cite{dum-tut-2017},  \cite{xie-2018}, \cite{malara-2018}). 

 For any homogenous ideal $I$ we can define the Castelnuovo-Mumford regularity as follows. Let $M$ be a finitely generated graded $S$-module. Let 
 \beqn
                             F_{\bullet}: 0 
 \longrightarrow F_{r}
 \longrightarrow F_{r-1}
 \cdots \longrightarrow F_1
 \longrightarrow F_0
 \longrightarrow 0
 \eeqn
be a minimal free resolution of $M$ where $F_{i} =\oplus_j S[-j]^{b_{ij}} $. Put $b_i(M) = \max\{ j | b_{ij} \not = 0\}$. Then  
$\reg(M) = \maxi\{ b_i(M) - i \}. $ 

Bounds on the Castelnuovo-Mumford regularity has been of interest. As the list is long  we state only a few results. 
In 1997,  Swanson  proved  that if $I$ is a homogenous ideal, then  there exists  an integer $r$ such that   
$reg(I^m) \leq mr$ for any m \cite{swanson-1997}. Later, it was proved that asymptotically  $\reg(I^m)$ is  a linear function of $m$ by  \cite{kodiyalam}, \cite{cutkosky-herzog-trung}.

The behaviour of  Castelnuovo-Mumford regularity of symbolic powers  is not easy to predict. From a 
  result of Cutkosky, Herzog and Trung, it follows that if $I$ is an ideal of points in a projective space and the symbolic Rees algebra $ R_s(I):\displaystyle{\bigoplus_{n \geq 0} } I^{(n)}$ is Noetherian, then  
  $\reg(I^{(n)})$ is a quasi-polynomial  (\cite[Theorem~4.3]{cutkosky-kurano}). Moreover, ${\displaystyle \lim_{n\rightarrow \infty}\left( \f{\reg(I^{(n)})}{n} \right)}$ exists and can even be irrational  \cite{cutkosky2}. For a nice survey article on Castelnuovo-Mumford regularity see \cite{chardin}.

Bocci and Harbourne showed  showed that if $I$ is  a  zero dimensional subscheme  in a projective space, then
   $\alpha(I) / {\gamma}(I) \leq \rho(I) \leq \reg (I) / {\gamma}(I)$ \cite[Theorem~1.2.1]{BH-2010}. Hence, if  
   $\alpha(I)=\reg(I)$, then $\rho(I)=\alpha(I)/{\gamma}(I)$.
   Later,   Harbourne and Huneke raised the following Conjecture: 
   
   \begin{conjecture}
    \cite[Conjecture~2.1]{harb-huneke-2013}
   Let $I$ be an ideal of fat points  in $S$ and $\m = (X_0, \ldots, X_d)$. Then   $I^{ (nd )  } \subseteq \m^{n(d-1)} I^n$ holds true for all $I$ and  $n$.
  \end{conjecture}
   
    In the same paper they  showed that the conjecture is true for fat point ideals arising as
  symbolic powers of radical ideals generated in a single degree in  $\Proj^2$.  
Recently, there has been a  renewed the interest on the Waldschmidt constant mainly due to the containment problem. In \cite{bocci-et-al-2016} the Waldschmidt constant for square free monomial ideals was computed. In fact they showed that if $\gamma(I)$ can be expressed as the value to a certain linear program arising from the structure of the associated primes of $I$. The Waldschmidt constant has also been computed for Stanley-Risner ideals \cite{bocci-franci-2016}. 
 
 The resurgence  and the Waldschmidt constant has been studied in a few cases: for certain general points in $\Proj^2$  \cite{BH0},  smooth subschemes  \cite{guardo}, fat linear subspaces 
 \cite{fatabbi-2015},  special point configurations \cite{duminicki} and monomial ideals \cite{bocco-waldschmidt-2016}.   
 
 We now briefly state our results on resurgence, Waldschmidt constant and Castelnuovo-Mumford regularity.
Putting  weights on  monomial curves ${\mathcal C}(\underline{a})$ in $\A^d$, we can consider them as weighted points in a weighted 
projective space $\Proj_{\kk}(\underline{a})$. Hence the bounds for resurgence in \cite{BH-2010} hold true. For 
$\underline{a} = (3,3+m,3+2m)$ these invariants have been computed in \cite{clare-shree}. For  $q  \geq 1$, and $\gcd(2q+1,m) = 1$ 
these invariants have been computed  $\underline{a}= (2q+1, 2q+1+m, 2q+1+2m)$, these invariants have been computed in 
\cite{clare}. In these cases the generators of the symbolic powers of $\p({\mathcal C}(\underline{a}))$ has been computed.
 In \cite{clare-mousumi}, for  monomial curves $\p({\overline{\mathcal C}}(\underline{a}))$ in $\Proj^3$, where $\underline{a}=   ( m, 2m, 2q+1+2m)$,  $q,m$ are positive integers and  
  $\gcd(2q+1,m)=1$,  these invariants have been computed.


\begin{thebibliography} {BB}
\bibitem
{BDHKKASS}
T.~Bauer,  S.~ Di Rocco, B.~ Harbourne, M.~ Kapustka, A.~ Knutsen, W.~ Syzdek and T.~Szemberg, 
{\em A primer on Seshadri constants.} Interactions of classical and numerical algebraic geometry, 33-70, 
Contemp. Math., {\bf 496}, Amer. Math. Soc., Providence, RI, 2009. 

\bibitem
{bocci-et-al-2016}
C.~Bocci,  S.~Cooper, E.~Guardo, Elena, B.~Harbourne, M.~Janssen,  U.~Nagel, Uwe, A.~Seceleanu, 
A.~ Van Tuyl, Adam and T. Vu, Thanh,
{\em The Waldschmidt constant for squarefree monomial ideals. }
J. Algebraic Combin. {\bf 44} (2016), no. 4, 875-904. 

\bibitem
{bocci-franci-2016}
C.~Bocci, Cristiano and B.~Franci, {\em Waldschmidt constants for Stanley-Reisner ideals of a class of simplicial complexes.} 
J. Algebra Appl. {\bf 15} (2016), no. 7, 1650137, 13 pp. 



\bibitem
{BH0}
C.~Bocci and  B.~Harbourne, {\em  The resurgence of ideals of points and the containment problem.}  Proc. Amer. Math. Soc. {\bf 138} (2010), no. 4, 1175-1190. 


 \bibitem
 {BH-2010} 
C.~Bocci and B.~Harbourne, {\em Comparing powers and symbolic powers of ideals}, J. Algebraic Geom. {\bf 19} (2010), no. 3, 399-417.


 \bibitem
 {bocco-waldschmidt-2016}
 C.~Bocci,  S.~Cooper,  E.~Guardo, B.~Harbourne, M.~Janssen,  U.~Nagel, A.~Seceleanu, A. ~Van Tuyl and Thanh Vu, {\em   The Waldschmidt constant for squarefree monomial ideals.} J. Algebraic Combin. {\bf 44} (2016), no. 4, 875-904.

\bibitem
{bor-1978}
M. Boraty{\' n}ski, {\em  A note on set-theoretic complete intersection ideals. }
J. Algebra {\bf 54} (1978), no. 1, 1-5. 


\bibitem
{bresinsky-1979}
H.~Bresinsky, H. Monomial space curves in $\A^3$ as set-theoretic complete intersections. Proc. Amer. Math. Soc. {\bf 75} (1979), no. 1, 23-24.

\bibitem
{bresinsky-a4-1979}
H.~Bresinsky, {\em Monomial Gorenstein curves in $\A^4$ as set-theoretic complete intersections.}  Manuscripta Math. {\bf 27} (1979), no. 4, 353-358. 

\bibitem
{chardin}
M.~Chardin,
{\em Some results and questions on Castelnuovo-Mumford regularity. Syzygies and Hilbert functions,} 1-40, 
Lect. Notes Pure Appl. Math., 254, Chapman \& Hall/CRC, Boca Raton, FL, 2007. 


 \bibitem
 {cowsik-1981}
  R.~C.~Cowsik,  {\em Symbolic powers and number of defining equations}. Algebra and its applications (New Delhi, 1981),  13-14, Lecture Notes in Pure and Appl. Math., 91, Dekker, New York, 1984. 

 \bibitem
{cowsik-nori} 
R.~C.~Cowsik and  M.~V.~Nori, {\em Affine curves in characteristic $\pp$ are set theoretic complete intersections}. Invent. Math. {\bf 45} (1978), no. 2, 111-114.

\bibitem
{cutkosky-1991}
S.~D.~Cutkosky, {\em Symbolic algebras of monomial primes.} J. Reine Angew. Math. {\bf 416} (1991), 71-89.

\bibitem
{cutkosky2}
S.~D.~Cutkosky,
{\em Irrational asymptotic behaviour of Castelnuovo-Mumford regularity.} 
J. Reine Angew. Math. {\bf 522} (2000), 93-103. 

\bibitem
{cutkosky-herzog-trung}
S~D.~Cutkosky, J. Herzog,  and  N. V. Trung, {\em  Asymptotic behaviour of the Castelnuovo-Mumford regularity.} Compositio Math. {\bf 118} (1999), no. 3, 243-261. 
 
 \bibitem
 {cutkosky-kurano}
 S. D. Cutkosky and K.~ Kurano,
{\em Asymptotic regularity of powers of ideals of points in a weighted projective plane. }
Kyoto J. Math. {\bf 51} (2011), no. 1, 25-45. 
 \bibitem
 {dao}
 H.~Dao,  A.~De Stefani, E.~Grifo, C.~Huneke and L.~N{\' u}{\~ n}ez-Betancourt, {\em  Symbolic powers of ideals. Singularities and foliations. geometry}, Topology and applications, 387-432, Springer Proc. Math. Stat., 222, Springer, Cham, 2018.

\bibitem
{clare}
C.~D'Cruz,
{\em Resurgence and Castelnuovo-Mumford regularity of certain monomial curves in $\A^3$.}  To appear in Acta Mathematica Vietnamica

\bibitem
{clare-shree}
C.~D'Cruz and S. Masuti, {\em Symbolic Blowup algebras and invariants of certain monomial curves in an affine space.} To appear in Comm. in Algebra. 

\bibitem
{clare-mousumi}
C.~D'Cruz and M.~Mandal, {\em Symbolic blowup algebras and invariants associated to certain monomial curves in $\Proj^3$.} To appear in Comm. in Algebra.

\bibitem
{duminicki}
M.~Dumnicki, B.~ Harbourne, U.~ Nagel, A.~ Seceleanu, T.~ Szemberg and  H.~Tutaj-Gasi{\' n}ska, 
{\em  Resurgences for ideals of special point configurations in $\Proj^{N}$ coming from hyperplane arrangements.} 
J.~Algebra {\bf 443} (2015), 383-394.

\bibitem
{dum-tut-2017}
M.~Dumnicki, and H.  Tutaj-Gasi{\' n}ska,
{\em A containment result in $\Proj^n$ and the Chudnovsky conjecture.} Proc. Amer. Math. Soc. {\bf 145} (2017), no. 9, 3689-3694. 

  \bibitem
 {ein-laz-smith-2001}
L.~Ein,  R.~Lazarsfeld, and K.~E.~ Smith,  {\em Uniform bounds and symbolic powers on smooth varieties}. Invent. Math. {\bf 144} (2001), no. 2, 241-252.

 \bibitem
 {eisenbud}
 D.~Eisenbud,  
 {\em Commutative algebra. With a view toward algebraic geometry.}  
 Graduate Texts in Mathematics, 150. Springer-Verlag, New York, 1995. 
 
 \bibitem
 {eisenbud-evans}
 D. Eisenbud and E. G. Evans  
 {\em Every algebraic set in $n$-space is the intersection of $n$ hypersurfaces.} 
 Invent. Math. {\bf 19} (1973), 107-112.



  \bibitem
  {fatabbi-2015}
 G.~ Fatabbi, Giuliana, B.~Harbourne and 
A.~Lorenzini, {\em Inductively computable unions of fat linear subspaces.} 
J. Pure Appl. Algebra {\bf 219} (2015), no. 12, 5413-5425.

\bibitem
{ferrand-1975}
 D. Ferrand, {\em Courbes gauches et fibr{\' e}s de rang 2. }
C. R. Acad. Sci. Paris S{\' e}r. A-B 281 (1975), no. 10, Aii, A345-A347.

\bibitem
{xie-2018}
L.~Fouli, Louiza,  P.~Mantero, Paolo and Y.~Xie,
{\em Chudnovsky's conjecture for very general points in $\Proj_{\kk}^N$ }
J. Algebra {\bf 498} (2018), 211-227.

 

\bibitem
{goto} 
 S.~Goto,  {\em The Cohen-Macaulay symbolic Rees algebras for curve singularities.} The Cohen-Macaulay and Gorenstein Rees algebras associated to filtrations. Mem. Amer. Math. Soc. {\bf 110} (1994), no.  526, 1-68. 
 
 \bibitem
{goto-nis-shim}
S.~Goto,  K. Nishida, and Y.~Shimoda,
{\em Topics on symbolic Rees algebras for space monomial curves.} Nagoya Math. J. {\bf 124} (1991), 99-132.
  
\bibitem
{goto-nis-shim-2}
S.~Goto,  K. Nishida,  Y.~Shimoda,
{\em The Gorensteinness of symbolic Rees algebras for space curves.}
J. Math. Soc. Japan {\bf 43} (1991), no. 3, 465-481.

\bibitem
{goto-nis-wat}
S.~Goto,  K. Nishida, K.~ Watanabe, {\em Non-Cohen-Macaulay symbolic blow-ups for space monomial curves and counterexamples to Cowsik's question.} Proc. Amer. Math. Soc. {\bf 120} (1994), no. 2, 383-392. 

 \bibitem
 {LNM-1092}
 S.~Greco and R.~Strano
 {\em Complete Intersections}
Lectures given at the 1st 1983 Session of the Centro Internationale Matematico Estivo (C.I.M.E.) held at Acireale (Catania), Italy, June 13-21, 1983.

\bibitem
{guardo}
E.~Guardo, B.~Harbourne and A.~ Van Tuyl, 
{\em Asymptotic resurgences for ideals of positive dimensional subschemes of projective space.} 
Adv. Math. {\bf 246} (2013), 114-127.

\bibitem
{harb-huneke-2013} 
B.~Harbourne and C.~Huneke,  {\em Are symbolic powers highly evolved?}
J. Ramanujan Math. Soc. {\bf 28A} (2013), 247-266. 

\bibitem
{hartshorne-LNM-156}
R.~Hartshorne, {\em Ample subvarieties of algebraic varieties.} Notes written in collaboration with C. Musili. Lecture Notes in Mathematics, Vol. 156 Springer-Verlag, Berlin-New York 1970 xiv+256 pp. 


\bibitem
{hartshorne-charp}
R.~Hartshorne, {\em Complete intersections in characteristic $\pp >0$}. Amer. J. Math. {\bf 101} (1979), no. 2, 380-383.

\bibitem
{herzog-1970}
J. Herzog: {\em Generators and relations of abelian semigroups and semigroup rings}. Manuscripta Math. 3 (1970) 175-193. 

 \bibitem
 {hochster-huneke-2002}
 M.~Hochster  and  C.~Huneke,
{\em Comparison of symbolic and ordinary powers of ideals}. 
Invent. Math. {\bf 147} (2002), no. 2, 349-369. 

\bibitem
{huneke-1982}
C.Huneke
On the finite generation of symbolic blow-ups.
Mathematische Zeitschrift {\bf 179} (1982) 465-472

\bibitem {huneke-1987} 
C.~Huneke,
{\em Hilbert functions and symbolic powers.}
Michigan Math. J. {\bf 34} (1987),   293-318. 


 \bibitem
 {schenzel-1992}
 G. Kn{\"o}del, P. Schenzel  and R. Zonsarow, {\em Explicit computations on symbolic powers of monomial curves in affine space. }
Comm. Algebra {\bf 20} (1992), no. 7, 2113-2126.

\bibitem
{kodiyalam}
V. Kodiyalam, {\em Asymptotic behaviour of Castelnuovo-Mumford regularity.} Proc. Amer. Math. Soc. {\bf 128} (2000), no. 2, 407-411. 

 \bibitem
 {kronecker}
 L.~Kronecker, {\em Zur Theorie der Abelschen Gleichungen.} J. Reine Angew. Math. {\bf 92}, 1-123 (1882).

\bibitem
{lyubeznik-survey}
 G.~Lyubeznik, 
{\em A survey of problems and results on the number of defining equations}. Commutative algebra (Berkeley, CA, 1987), 375-390, Math. Sci. Res. Inst. Publ., 15, Springer, New York, 1989. 
 
 \bibitem
 {lyubeznik-1992}
 G.~Lyubeznik, {\em The number of defining equations of affine algebraic sets. }
Amer. J. Math. {\bf 114} (1992), no. 2, 413-463. 

\bibitem
{malara-2018}
G.~Malara, T.~Szemberg, and J.~Szpond, 
{\em On a conjecture of Demailly and new bounds on Waldschmidt constants in $\Proj^N$.}  
J. Number Theory {\bf 189} (2018), 211-219. 

 \bibitem
 {moh-1982}
T.~ T.~Moh,
{\em A result on the set-theoretic complete intersection problem. }
Proc. Amer. Math. Soc. {\bf 86} (1982), no. 1, 19-20. 

 \bibitem
 {moh-1985}
 T.~T.~Moh
 {\em Set-theoretic complete intersections. }
Proc. Amer. Math. Soc. {\bf 94} (1985), no. 2, 217-220. 
 
 \bibitem
 {kumar-1978}
 N. Mohan Kumar, 
{\em On two conjectures about polynomial rings. }
Invent. Math. {\bf 46} (1978), no. 3, 225-236.  
 
 
 \bibitem{morales-1991}
 M.~Morales, {\em Noetherian symbolic blow-ups.} J. Algebra {\bf 140} (1991), no. 1, 12-25.


 \bibitem
 {nagata}
M.~Nagata. {Local rings} Interscience Tracts in Pure and Applied Mathematics, No. {\bf 13} Interscience Publishers a division of John Wiley \& Sons. New York-London 1962 xiii+234 pp.
 
 \bibitem
 {nagata-1965}
 M.~Nagata
 {\em Lectures on the fourteenth problem of Hilbert.}
  Tata Institute of Fundamental Research, Bombay 1965 ii+78+iii pp. 
  
 \bibitem
 {patil-1990}
 D. Patil, {\em Certain monomial curves are set theoretic complete intersection.} Manuscripta Math. {\bf 68}, (1990)
 399-404.
 
 \bibitem
 {rees-1958}
  D.~Rees, D. {\em On a problem of Zariski.} Illinois J. Math. {\bf 2} 1958 145-149. 
    
 \bibitem
 {robbiano-valla-1983}
 Robbiano L., Valla G. {\em Some curves in $\Proj^3$ are set-theoretic complete intersections.}  Algebraic geometry-open problems (Ravello, 1982), 391-399, Lecture Notes in Math., 997, Springer, Berlin-New York, 1983.

  \bibitem
  {roberts-1985}
  P.~Roberts,   {\em  A prime ideal in a polynomial ring whose symbolic blow-up is not Noetherian. } Proc. Amer. Math. Soc. {\bf 94} (1985), no. 4, 589-592. 


 \bibitem
 {sahin-2009-pn} {\c{S}}ahin, Mesut,  {\em Producing set-theoretic complete intersection monomial curves in $\Proj^n$}. Proc. Amer. Math. Soc. {\bf 137} (2009), no. 4, 1223-1233.
 
 \bibitem
 {schenzel-1988}
 P. Schenzel, Examples of Noetherian symbolic blow-up rings, Rev. Rom. Math Pures et appl. {\bf 33} (1988) 375-383.
  
\bibitem
{hema-1991}
H.~Srinivasan,  {\em On finite generation of symbolic algebras of monomial primes. }
Comm. Algebra {\bf 19} (1991), no. 9, 2557-2564. 
 
 \bibitem
 {storch}
 U.~Storch, {\em Bemerkung zu einem Satz} von M. Kneser.  Arch. Math. (Basel) {\bf 23} (1972), 403-404.

\bibitem
{stuckard-vogel-1982}
J.~St{\"u}ckrad and W.~Vogel,
{\em On the number of equations defining an algebraic set of zeros in
$n$-space.} Seminar D. Eisenbud/B. Singh/W. Vogel, Vol. 2, pp. 88-107, 
Teubner-Texte zur Math., 48, Teubner, Leipzig, 1982. 

\bibitem
{swanson-1997}
I. Swanson, {\em Powers of ideals. Primary decompositions, Artin-Rees lemma and regularity.} Math. Ann. {\bf 307} (1997), no. 2, 299-313.

\bibitem
{swanson-2000}
I. Swanson, {\em Linear equivalence of ideal topologies.}  Math. Z. {\bf 234} (2000), no. 4, 755-775. 

\bibitem
{szpiro-1979}
L.~Szpiro, "Lectures on equations defining space curves," 
Notes by N. Mohan Kumar. Tata Institute of Fundamental Research, Bombay; by Springer-Verlag, Berlin-New York, (1979).

 \bibitem
 {valla-1980}
 G.~Valla, 
{\em On determinantal ideals which are set-theoretic complete intersections. }
Compositio Math. {\bf 42} (1980/81), no. 1, 3-11. 

 \bibitem
 {waldschmidt}
M.~ Waldschmidt,  {\em Propri{\' e}t{\' e}s arithm{\' e}tiques de fonctions de plusieurs variables. II.}
 (French)  S{\' e}minaire Pierre Lelong (Analyse) ann{\' e}e 1975/76, pp. 108-135. 
Lecture Notes in Math., Vol. 578, Springer, Berlin, 1977. 

 \bibitem
 {zariski}
 O~Zariski. {\em A fundamental lemma from the theory of holomorphic functions on an algebraic variety. }
Ann. Mat. Pura Appl. (4) {\bf 29} (1949), 187-198. 
 
   \end{thebibliography}
\end{document}